\documentclass[a4paper, reqno, 11pt]{amsart}

\usepackage[a4paper,margin=25mm]{geometry}
\usepackage{amsmath,amssymb,amsthm,mathtools}
\usepackage{microtype}
\usepackage[hidelinks]{hyperref}
\usepackage{xcolor}

\newtheorem{theorem}{Theorem}
\newtheorem{lemma}{Lemma}
\newtheorem{remark}{Remark}

\newtheorem{definition}{Definition}

\newcommand{\R}{\mathbb R}
\newcommand{\T}{\mathbb T}
\newcommand{\nuav}[1]{\langle\nu,#1\rangle}
\newcommand{\dd}{\,\mathrm d}
\newcommand{\Om}{\Omega}
\newcommand{\QT}{Q_T}
\newcommand{\del}{\partial}

\DeclareMathOperator{\diver}{div}
\DeclareMathOperator{\curl}{curl}

\begin{document}

\title[Adiabatic Thermoelasticity]{An energy-conserving variational scheme \\ for adiabatic
thermoelasticity  \\ with convex stored energy}

\author[K. Koumatos]{Konstantinos Koumatos}
\address[Konstantinos Koumatos]{\newline
Department of Mathematics
\newline
University of Sussex
\newline
Pevensey 2 Building
\newline
Falmer, Brighton, BN1 9QH, UK
}
\email{k.koumatos@sussex.ac.uk}

\author[A.E. Tzavaras]{Athanasios E. Tzavaras}
\address[Athanasios E. Tzavaras]{
\newline
Computer, Electrical and Mathematical Sciences and Engineering Division 
\newline
King Abdullah University of Science and Technology (KAUST)
\newline 
Thuwal 23955-6900,  Saudi Arabia
}
\address{
and 
\newline
Institute of Applied and Computational Mathematics
\newline
Foundation of Research and Technology - Hellas (FORTH)
\newline 
Heraklion, Greece
}
\email{athanasios.tzavaras@kaust.edu.sa}

\author[A. Vikelis]{Andreas Vikelis}
\address[Andreas Vikelis]{\newline 
Department of Mathematics
Faculty of Mathematics
University of Vienna
Oskar-Morgenstern-Platz 1
1090 Vienna, Austria}
\email{andreas.panagiotis.vikelis@univie.ac.at}

\date{}

\vspace{-2em}

\begin{abstract}
We propose an approximation scheme for the equations of adiabatic thermoelasticity in several space dimensions
that is motivated by the  thermodynamical structure of the system. The approximating scheme conserves
 the total energy  and captures entropy production consistent with the Clausius-Duhem inequality.
As an application, we prove existence of dissipative measure-valued solutions to the system of 
adiabatic thermoelasticity on the torus with convex energy and prescribed nonnegative heat supply. {The constructed dissipative
measure-valued solutions satisfy an averaged form of the energy conservation.}
\end{abstract}

\maketitle

\tableofcontents

\section{Introduction}

The system of adiabatic thermoelasticity in Lagrangian coordinates 
consists of the balance laws of momentum and energy
\begin{equation}\label{smooth-system}
\begin{aligned}
 \partial_t F &= \nabla v \, ,
 \\
 \partial_t v &= \diver S(F,\eta) \, ,
 \\
 \partial_t \Big ( \tfrac{1}{2} |v|^2 + e(F,\eta) \Big )
 &= \diver\big( v^T S(F,\eta) \big)+r.
 \end{aligned}
\end{equation}
Here $F=\nabla y$, $v=\partial_t y$, $\eta$ is the specific entropy, and $F$ 
satisfies the constraint $\curl F = 0$. The quantity $r$ is called radiative heat supply
and we omit body forces for simplicity.
The system describes the evolution of
a thermomechanical process $(v, F, \eta)(t,x)$ taking values
in $\mathcal Z=\R^d\times\R^{d\times d}\times[0,\infty)$.
The constitutive relations
\begin{equation}
S(F,\eta) = \frac{\del e}{\del F} (F,\eta) \quad \mbox{and} \quad  \theta(F,\eta) = \frac{\del e}{\del \eta} (F,\eta) >0 
\end{equation}
make the resulting system consistent with the Clausius-Duhem inequality,
which is a non-equilibrium version of the second law of thermodynamics \cite{TN1992}.
As a result, by design, smooth solutions of \eqref{smooth-system} satisfy the entropy production identity 
\begin{equation}\label{entropy-cons}
\partial_t\eta= \frac{r}{\theta(F,\eta)}  \, .
\end{equation}

{
System \eqref{smooth-system} is one of the standard paradigms in the theory of hyperbolic conservation laws,
see \cite{b-Dafermos2026}. Hyperbolicity of that system requires the strict Legendre--Hadamard condition in $F$, 
namely 
$$
D^2_{FF}e(F,\eta)[a\otimes n,a\otimes n]>0 \quad \mbox{for $a,n\ne0$}
$$
together with $e_\eta \ne 0$. Instead, we work with a relaxed set of hypotheses
\begin{equation}\label{assumption}
D^2_{FF}e(F,\eta) > 0  \quad \mbox{ and} \quad \theta = e_\eta (F, \eta) > 0 \, ,
\tag {A}
\end{equation}
which are adequate in explaining the methodology, and will relax to more general convexity
hypotheses in a future work.
}
The one-dimensional variant of \eqref{smooth-system}, 
\begin{equation}\label{1dsystem}
\begin{aligned}
 \partial_t u &= \del_x v \, ,
 \\
 \partial_t v &= \del_x \sigma (u,\eta) \, ,
 \\
 \partial_t \Big ( \tfrac{1}{2} |v|^2 + e(u,\eta) \Big )
 &= \diver\big(  \sigma (u,\eta) v \big)+r \ ,
 \end{aligned}
\end{equation}
is a model for adiabatic gas dynamics in Langrangian coordinates.
The thermomechanical process is described by $(u,v,\eta)(t,x)$, taking values in $(0,\infty)\times \R \times[0,\infty)$,
but  $u$ is now interpreted as the specific volume and is required to be positive $u >0$.

The second law of thermodynamics plays a key role in the theory of thermodynamic processes near equilibrium
with important implications in gas dynamics and chemical physics, see Callen \cite{b-Callen1985}.
In the theoretical mechanics literature a non-equilibrium version of the second law, often
referred to as the Clausius-Duhem inequality, has been used extensively as a
design principle requiring that models of continuum physics are
automatically endowed with additional structure, see the treatise \cite{TN1992}
and the works \cite{CN1963,CM1964} for implications on the theory of thermoviscoelasticity. 
The interested reader may consult the surveys \cite{Silhavy1986}, \cite{Serrin1996}, \cite{Evans} for a
wealth of information on the interface of thermodynamics and partial differential equations.

This fruitful interface was instrumental in the theory of hyperbolic conservation laws.
The connection motivated Lax \cite{Lax1971}, \cite{FL1971} to introduce the notion
of entropy solution as a mechanism for the selection of admissible non-smooth processes and 
to study its implications on the stability of numerical schemes.
A lot of work ensued on both the role of the second law as a stability mechanism, see Dafermos \cite{Dafermos1979},
and its paramount implications on the design and stability properties of numerical schemes for systems
of conservation laws, {\it e.g.} Tadmor \cite{Tadmor2003}. The reader may consult the treatise \cite{b-Dafermos2026} 
for multiple aspects of this beneficial interface between the field of natural philosophy and the 
subject of conservation laws.

Our interest here is in the interface between the calculus of variations and aspects of the dynamical theory of
systems of conservation laws. It is well known that several special systems of mechanics and mathematical physics
- mostly isothermal - can be derived from the principle of least action and often have a Hamiltonian structure. 
An ongoing question is whether such systems can be approached via techniques from the calculus of variations while respecting
the entropy structure for weak solutions. 
Various works have been devoted to isothermal elastodynamics in Lagrangian coordinates
\cite{DST2001, DST2012, KS2019}, and a parallel line of work exists on the construction of weak solutions 
for isentropic gas dynamics in Eulerian coordinates via techniques of the calculus of variations,
\cite{GW2009,Westdickenberg2017}.{ Entropy weak solutions are produced in one-space dimension \cite{DST2000},
but  only measure-valued solutions in several space dimensions, due to the well known theoretical challenges for hyperbolic 
conservation laws.}

{
Later, attention turned to nonisothermal equations, with the aim of extending the
aforementioned works to the full thermomechanical system; see \cite{CGT2019,CGT2020,GKV2022} for
the system of thermoelasticity in Lagrangian coordinates and  \cite{CSW2019} for the compressible Euler equations in Eulerian coordinates.
These works have the appealing feature that they offer construction schemes for the relevant systems that may have
implications in future designs of stable algorithms.
A drawback of the above works on non-isothermal systems is that the equation of energy is often deduced as an energy inequality 
rather than as an energy identity. It may be that in some processes (for instance turbulence) energy inequalities may play an important role.
Moreover, the study of non-equilibrium processes where the first law is viewed as an inequality has received some attention in the theoretical mechanics
literature, see \cite{FS1975,Serrin1996} for the role of an inequality in the first law and the discussion in \cite{CGT2020} on implications on
the mathematical theory of thermomechanical systems. Nevertheless, in many thermomechanical processes the energy is conserved and the entropy is produced
and this "more conventional" case has not been adequately analyzed.
}

The objective of this work is to address how to rectify this drawback, and it will turn out that thermodynamics will play a key role
in bridging the gap.  There are two origins of this difficulty. One is of a technical nature:
the energy may develop a concentration defect because the available bounds do not ensure uniform integrability of the energy densities.
A second difficulty is that approximation schemes may contain errors that produce defect measures in the limit.
The objective is, by exploiting the underlying thermodynamical structure, to propose an approximating scheme that conserves the energy
and is consistent with entropy production. We use the scheme to construct a dissipative measure-valued solution for the system of adiabatic thermoelasticity
\eqref{smooth-system} (see Definition \ref{defDMV} for the precise definition of a dissipative measure-valued solution).
Its energy will have an averaged part connected to a Young measure and a defect part but the total energy will satisfy
a space-integrated energy identity, see Theorem \ref{thm} for the statement. The averaged entropy production satisfies an entropy inequality;
this measure-valued entropy inequality includes the averaged reciprocal temperature.

{
The present work is the first in a series that will derive appropriate minimizing schemes consistent with thermodynamics and energy conservation for 
the system \eqref{smooth-system}. In future work we will address the full class of polyconvex or quasiconvex energies. In the current study, in order to
simplify the exposition and give a complete theory,  we study the simplest case \eqref{assumption} under
the additional restriction $r \ge 0$. 
}
The proposed approximation scheme is motivated by the combined action of the first and second laws of thermodynamics in the approximation process. 
Recall that  in equilibrium thermodynamics a central role is played 
by the study of elementary quasistatic thermodynamic processes (such as isothermal, isobaric or adiabatic processes). Likewise, in the present context which
is far from equilibrium, we may consider the effect of various non-equilibrium elementary processes. In the scheme we propose we will alternate between
two such elementary processes: (i) one with fixed entropy that has a variational implementation;
(ii) a second with fixed mechanical characteristics that updates the entropy in such a way that the total energy is conserved. We note that entropy updates are also used in  \cite{del_pino_1,del_pino_2} for developing convex optimization algorithms for the case of gas dynamics.

The outline of this work is as follows:
After some preparation in Section \ref{sec:prelim}, we develop the approximation scheme in Section \ref{sec:constrvarscheme} and show its solvability and basic estimates. 
Finally, in Section \ref{sec:construction} we construct the dissipative measure-valued solution and prove Theorem \ref{thm}.

\section{Preliminaries}\label{sec:prelim}

\subsection{The constitutive theory of thermoviscoelasticity}

It is instructive to first review the constitutive theory of  thermoviscoelasticity.
These equations are derived in \cite{CN1963,CM1964} by asking what restrictions need to be imposed on general constitutive theories
of thermoviscoelasticity so that they become consistent with the Clausius-Duhem inequality. The emerging system  has the form:
\begin{equation}\label{system-ve}
\begin{aligned}
 \partial_t F &= \nabla v \, ,
 \\
 \partial_t v &= \diver \big ( S(F,\eta) + Z(F, \eta, \nabla v, \nabla \theta) \big ) \, ,
 \\
 \partial_t \Big ( \tfrac{1}{2} |v|^2 + e(F,\eta) \Big )
 &= \diver\left( v^T \big ( S(F,\eta) + Z(F, \eta, \nabla v, \nabla \theta) \big ) \right )+ \diver Q(F,\eta,\nabla v, \nabla \theta) + r.
 \end{aligned}
\end{equation}
The quantities $F$ and $v$ take values in $\R^{d\times d}$ and $\R^d$, respectively, while the scalar-valued entropy
is assumed bounded from below; for simplicity we take $\eta\geq0$.
The total stress $T \in\R^{d\times d}$ is decomposed into an elastic part $S(F,\eta)$ and a viscoelastic part $Z(F, \eta, \nabla v, \nabla \theta)$
with $Z(F,\eta, 0, \nabla \theta) = 0$, according to the formula
$$
T = S(F,\eta) + Z(F, \eta, \nabla v, \nabla \theta) \, ,
$$
where the elastic part of the stress $S(F,\eta)$ and the temperature $\theta (F, \eta)$  are derived from a potential function which is the internal
energy  $e(F,\eta)$: 
\begin{equation}\label{maxwell}
 S(F,\eta)=D_Fe(F,\eta),\qquad
 \theta(F,\eta)=\partial_\eta e(F,\eta).
\end{equation}
Finally it is assumed that $Z$ and the heat flux $Q$ satisfy the inequalities
\begin{equation}\label{dissip}
Z(F, \eta, \nabla v, \nabla \theta)  : \nabla v \ge 0 \quad \mbox{and} \quad Q(F,\eta,\nabla v, \nabla \theta) \cdot \nabla \theta \ge 0
\end{equation}
The reader may keep in mind the special selections $Z = \mu (F,\eta) \nabla v$ of Newtonian viscosity
and $Q(F,\eta,\nabla v, \nabla \theta) = \kappa (F, \eta) \nabla \theta$ of Fourier heat conduction, in which case
the requirement \eqref{dissip} reduces to the conditions  $\mu \ge 0$ and $\kappa \ge 0$ for the viscosity and heat conductivity.

The following calculation monitors the evolution of the entropy. For smooth solutions, subtracting from the total energy \eqref{system-ve}$_3$
the kinetic energy (formed by multiplying $v$ with \eqref{system-ve}$_2$) gives
\begin{equation}\label{energydiss}
\del_t e (F,\eta) = (S+Z) : \nabla v + \diver Q + r
\end{equation}
For smooth processes, a simplification using \eqref{maxwell} leads to
\begin{equation}\label{entropy}
\begin{aligned}
\del_t \eta - \diver \frac{Q}{\theta} &= \frac{1}{\theta} Z : \nabla v + \frac{1}{\theta^2} Q \cdot \nabla \theta + \frac{r}{\theta}
\\
&= D + \frac{r}{\theta} \ge \frac{r}{\theta}
\end{aligned}
\end{equation}
where $D$ stands for the effect of the viscous and thermal dissipation which is nonnegative by \eqref{dissip}.

On the other hand the total energy (kinetic plus potential) is conserved and if for example $\Om=\T^d$ is the torus we see
that
\begin{equation}\label{energy}
E(t) := \int_\Om \tfrac{1}{2} |v|^2 + e(F,\eta) dx = E(0) + \int_0^t \int_\Om r(s,x)  dx ds
\end{equation}

Consider now the system of adiabatic thermoelasticity \eqref{smooth-system}, the topic of study in this work.
Motivated by viewing thermoelasticity as a limiting case of the theory of thermoviscoelasticity, system \eqref{smooth-system} may be regarded as emerging from a limiting process of \eqref{system-ve}, whereby solutions of adiabatic thermoelasticity are
obtained  by sending the viscosity and heat conduction to zero, $Z \to 0$, $Q\to 0$. 
One observes that, in the limit, we will still obtain the conservation of energy equation \eqref{energy} while the entropy identity \eqref{entropy} becomes 
an entropy production inequality
\begin{equation}\label{entropy-prod}
\partial_t\eta \ge  \frac{r}{\theta(F,\eta)}  \, .
\end{equation}
Our objective in Section \ref{sec:constrvarscheme} is to propose an alternative approximation scheme 
that is variational in conception. When pursuing this scheme we will pay close attention to 
preserve the properties \eqref{energy} and \eqref{entropy-prod} along the process of iteration.

\subsection{Hypotheses}

We begin with listing the framework and the hypotheses employed.
Let $\Om=\T^d$ be the $d$-dimensional torus, equipped with its Lebesgue measure, and let
$T>0$ and $\QT=(0,T)\times\Om$. All spatial test functions will be periodic. 
We use the convention
\[
 (\nabla v)_{i\alpha}=\partial_\alpha v_i,
 \qquad
 {\color{black}(\diver F)_i=\sum_\alpha\partial_\alpha F_{i\alpha},\qquad (v\cdot S)_\alpha=(v^TS)_\alpha=\sum_i v_iS_{i\alpha}.}
\]

\vspace{0.2cm}

\noindent
For the energy function $e(F,\eta)$ we assume it satisfies the hypotheses:
\begin{enumerate}
\item[(H1)]
{ $e\in C^1 (\R^{d\times d}\times[0,\infty))$}

\item[(H2)]
There are $c_0,c_1,c_2>0$ such that
\[
 c_0(|F|^2+\eta^2)-c_1
 \leq e(F,\eta)
 \leq c_2(1+|F|^2+\eta^2).
\]

\item[(H3)]
There is $\kappa>0$ such that, for every
$F,G\in\R^{d\times d}$ and $\eta\geq0$,
\[
 e(G,\eta)
 \geq e(F,\eta)+S(F,\eta):(G-F)
       +\frac{\kappa}{2}|G-F|^2.
\]
Thus the energy is uniformly convex in the deformation
gradient.

\item[(H4)]
There is $\theta_0>0$ such that
$\theta(F,\eta)\geq\theta_0$ for all $F$ and $\eta\geq0$.

\item[(H5)]
There is $C_S>0$ such that
\[
 |S(F,\eta)|\leq C_S(1+|F|+\eta).
\]

\item[(H6)]
There is $C_\theta>0$ such that
\[
 \theta(F,\eta)\leq C_\theta(1+|F|+\eta).
\]
\end{enumerate}

\subsection{Total and relative energy}

Following the total-energy notation of \cite{CGT2020},
let $U=(v,F,\eta)$ and define
\begin{equation}\label{Idefinition}
 I(U)=I(v,F,\eta)=\frac12|v|^2+e(F,\eta).
\end{equation}
For $U_*=(v_*,F_*,\eta_*)$, define the Bregman/Taylor remainder
\begin{align}\label{relativeI}
 I(U_*\mid U)
 &=I(U_*)-I(U)-DI(U)\cdot(U_*-U)\\
 &=\frac12|v_*-v|^2+e(F_*,\eta_*)-e(F,\eta)
 -S(F,\eta):(F_*-F)
 -\theta(F,\eta)(\eta_*-\eta).
 \notag
\end{align}
Only equal-entropy comparisons are needed for the
construction. For $\eta_*=\eta=a$, (H3) gives
\begin{equation}\label{relativepositive}
 I(v_*,F_*,a\mid v,F,a)
 \geq\frac12|v_*-v|^2+\frac{\kappa}{2}|F_*-F|^2.
\end{equation}
No positivity is asserted for \eqref{relativeI} at
different entropy values without joint convexity of $e$.
This distinction is important: the existence hypotheses
here require uniform convexity only in $F$.

Choose $C_*=c_1+1$ and set $\widetilde I=I+C_*$.
By (H2), for suitable constants $c,C>0$,
\begin{equation}\label{shiftgrowth}
 c(1+|v|^2+|F|^2+\eta^2)
 \leq\widetilde I(v,F,\eta)
 \leq C(1+|v|^2+|F|^2+\eta^2).
\end{equation}
That is, shifting the energy leaves all balance laws and
relative-energy remainders unchanged.

\subsection{Dissipative measure-valued solutions}

Let us begin by observing that the periodic motion used here imposes $\int_\Omega F^0\,dx=0$.
For a prescribed nonzero mean gradient, one may instead write
$y=F_{\mathrm{av}}x+u$ with periodic $u$; this extension is not made here.
Let
\[
 y^0\in W^{1,2}(\Om;\R^d),\qquad
 F^0=\nabla y^0,\qquad
 v^0\in L^2(\Om;\R^d),\qquad
 \eta^0\in L^2(\Om),\quad\eta^0\geq0.
\]
Put
\[
 E_0=\int_\Om I(v^0,F^0,\eta^0)\dd x,
 \qquad
 \nu_{0,x}=\delta_{(v^0(x),F^0(x),\eta^0(x))},
 \qquad
 \gamma_0=0,
\]
and define
\[
 \mathcal E(t)
 =E_0+\int_0^t\int_\Om r(s,x)\dd x\dd s.
\]
Write
\[
 \mathcal Z=\R^d\times\R^{d\times d}\times[0,\infty).
\]

\begin{definition}\label{defDMV}
A dissipative measure-valued solution with these initial
data and prescribed heat supply $r$ consists of a motion
$y$, a probability Young measure $\nu$ on $\mathcal Z$,
and a nonnegative finite Radon measure $\gamma$ on $Q_T$,
with the following properties.
\begin{enumerate}
\item
The moments of $\nu$ satisfy
\[
 \operatorname*{ess\,sup}_{t\in(0,T)}
 \int_\Om\nuav{|F|^2+|v|^2+\eta^2}\dd x<\infty.
\]
The barycentres
$\bar F=\nuav{F}$, $\bar v=\nuav{v}$ and
$\bar\eta=\nuav{\eta}$ satisfy
\[
 y\in W^{1,\infty}(0,T;L^2)
       \cap L^\infty(0,T;W^{1,2}),
\]
\[
 \nabla y=\bar F,\qquad
 \partial_t y=\bar v,\qquad
 y(0)=y^0.
\]

\item
There is a generating sequence
$(v_k,\nabla y_k,\eta_k)$ for $\nu$, with
$v_k=\partial_t y_k$, $\eta_k\geq0$, and uniform bounds
\begin{align*}
 &\|y_k\|_{L^\infty W^{1,2}}
 +\|v_k\|_{L^\infty L^2}
 +\|\eta_k\|_{L^\infty L^2}\leq C,\\
 &\|\partial_t\nabla y_k\|_{L^\infty H^{-1}}\leq C.
\end{align*}
For this sequence, with $\widetilde I=I+C_*\geq0$,
the energy measures satisfy
\[
 \widetilde I(v_k,\nabla y_k,\eta_k)\dd x\dd t
 \overset{*}{\rightharpoonup}
 \nuav{\widetilde I}\dd x\dd t+\gamma.
\]

\item
The averaged equations
\[
 \partial_t\bar F=\nabla\bar v,\qquad
 \partial_t\bar v=\diver\nuav{S(F,\eta)},\qquad
 \partial_t\bar\eta
 \geq r\nuav{1/\theta(F,\eta)}
\]
hold in distributions, with the prescribed initial terms
in their weak formulations. Explicitly, the mechanical
weak forms are \eqref{weakF}--\eqref{weakv} below, and
\[
 -\int_{Q_T}\bar\eta\partial_t\zeta\dd x\dd t
 -\int_\Om\eta^0\zeta(0)\dd x
 \geq
 \int_{Q_T}r\nuav{1/\theta(F,\eta)}\zeta\dd x\dd t
\]
for every nonnegative smooth periodic $\zeta$ compactly
supported in $[0,T)\times\Om$.

\item
For every nonnegative $\varphi\in C_c^1([0,T))$,
\begin{equation}\label{DMVenergy}
 E_0\varphi(0)
 +\int_{Q_T}\varphi'(t)\nuav{I}\dd x\dd t
 +\int_{Q_T}\varphi'(t)\dd\gamma
 +\int_{Q_T}r\varphi\dd x\dd t=0.
\end{equation}
\end{enumerate}
\end{definition}

\section{A constrained variational scheme}\label{sec:constrvarscheme}

\noindent We construct the discrete solution in two steps. For simplicity, we will assume that $r (t,x) \ge 0$.
First, we
update the velocity and the deformation gradient by solving a constrained
minimization problem while keeping the entropy at its previous value.
The mechanical update decreases the total energy by a nonnegative
quantity consisting of a kinetic increment and an elastic Bregman
remainder. In the second step, we restore this energy through a
pointwise increase of the entropy. The prescribed nonnegative heat
supply enters the same correction. Consequently, the completed update
satisfies an exact spatially integrated total-energy balance and has
nonnegative entropy production.

\vspace{0.2cm}

\noindent We retain the notation $\Om=\mathbb T^d$ and the
constitutive assumptions of the preliminaries.
All functions and test functions are periodic in the spatial
variable. Unless otherwise specified, the spaces $L^p$ and
$W^{1,2}$ are taken over $\Om$, with scalar, vector,
or matrix values as appropriate. We write
\[
H^{-1}(\Om;\mathbb R^d)
=
\bigl(W^{1,2}(\Om;\mathbb R^d)\bigr)^*.
\]
{Spatial} distributional identities are tested against smooth
periodic functions of the corresponding dimension. Pointwise
identities and inequalities between measurable functions are
understood almost everywhere.

\vspace{0.2cm}

\noindent Recall that
\[
S(F,\eta)=\partial_F e(F,\eta),
\qquad
\theta(F,\eta)=\partial_\eta e(F,\eta).
\]
In particular, $e$ is continuously differentiable, with the
derivative in $\eta$ understood from the right at $\eta=0$, and
\begin{align}
c_0\bigl(|F|^2+\eta^2\bigr)-c_1
&\leq e(F,\eta)
\leq c_2\bigl(1+|F|^2+\eta^2\bigr),
\label{eq:energy-growth}\\
|S(F,\eta)|+\theta(F,\eta)
&\leq C(1+|F|+\eta),
\label{eq:derivative-growth}\\
\theta(F,\eta)&\geq\theta_0>0.
\label{eq:temperature-lower}
\end{align}
Here and below, $\eta\geq0$. The uniform convexity assumption
in the deformation gradient reads
\begin{equation}
e(G,\eta)
\geq e(F,\eta)+S(F,\eta):(G-F)
+\frac{\kappa}{2}|G-F|^2,
\qquad \kappa>0.
\label{eq:uniform-convexity}
\end{equation}
No convexity in $\eta$ is needed in the entropy correction below.

\vspace{0.2cm}

\noindent\underline{\textbf{The mechanical minimization.}}
Fix $T>0$, let $h=T/N$, and set $t_j=jh$.
For the prescribed heat supply
\[
r\in L^1\bigl((0,T)\times\Om\bigr),
\qquad r\geq0,
\]
define
\begin{equation}
r^j(x)=\frac1h\int_{t_{j-1}}^{t_j}r(t,x)\,dt.
\label{eq:heat-average}
\end{equation}
Then $r^j\in L^1(\Om)$ and $r^j\geq0$ almost everywhere.

\vspace{0.2cm}

\noindent Suppose that the previous state satisfies
\[
\begin{gathered}
v^{j-1}\in L^2(\Om;\mathbb R^d),
\qquad
F^{j-1}\in L^2(\Om;\mathbb R^{d\times d}),\\
\eta^{j-1}\in L^2(\Om),
\qquad
\eta^{j-1}\geq0.
\end{gathered}
\]
The initial state $(v^0,F^0,\eta^0)$ is assumed to satisfy
these same $L^2$ conditions, with $\eta^0\geq0$.
The construction below shows that these conditions are
preserved at each step. The spatial Sobolev regularity
$v^j\in W^{1,2}(\Om;\mathbb R^d)$ is obtained
for $j\geq1$; no such regularity is required of $v^0$.

\vspace{0.2cm}

\noindent The mechanical unknowns are a pair $(w,G)$ subject
to the discrete kinematic constraint. Define
\begin{equation}
\mathcal A_j
=
\left\{
\begin{aligned}
&(w,G)\in W^{1,2}(\Om;\mathbb R^d)
\times L^2(\Om;\mathbb R^{d\times d}):
&
\frac{G-F^{j-1}}h=\nabla w
\end{aligned}
\right\},
\label{eq:admissible-class}
\end{equation}
where the constraint is affine and understood in the sense
of distributions.

\vspace{0.2cm}

\noindent More explicitly, the distributional constraint means that
\[
\int_\Om(G-F^{j-1}):\Phi\,dx
=
-h\int_\Om w\cdot\operatorname{div}\Phi\,dx
\]
for every
$\Phi\in C^\infty(\Om;\mathbb R^{d\times d})$,
where
$(\operatorname{div}\Phi)_i
=\sum_\alpha\partial_\alpha\Phi_{i\alpha}$.
Since $G-F^{j-1}$ and $\nabla w$ both belong to $L^2$,
this constraint is also an equality in
$L^2(\Om;\mathbb R^{d\times d})$ and hence
holds almost everywhere. By density, its weak formulation
also holds for
$\Phi\in W^{1,2}(\Om;\mathbb R^{d\times d})$.

\vspace{0.2cm}

\noindent We choose $(v^j,F^j)$ to minimize
\begin{equation}
\mathcal J_j(w,G)
=
\int_\Om
\left(
\frac12|w-v^{j-1}|^2+e(G,\eta^{j-1})
\right)\,dx
\label{eq:functional}
\end{equation}
over $\mathcal A_j$. Thus the minimization problem is
\begin{equation}
(v^j,F^j)\in
\operatorname*{arg\,min}_{(v,F)\in\mathcal A_j}
\mathcal J_j(v,F),
\label{eq:minimization}
\end{equation}
where, in particular,
\begin{equation}
F^j-F^{j-1}=h\nabla v^j
\label{eq:kinematic-update}
\end{equation}
in the sense of distributions. The entropy $\eta^{j-1}$ is a
prescribed coefficient in this problem.
Keeping it fixed ensures that the variation of the internal
energy produces precisely the stress evaluated at $\eta^{j-1}$, i.e.
$S(F^j,\eta^{j-1})$. The entropy is updated only after the
mechanical dissipation has been identified.

\vspace{0.2cm}

\begin{lemma}
Fix $j\in\mathbb N$. The minimization problem
\eqref{eq:minimization} is well defined and admits a unique
minimizer.
\end{lemma}

\begin{proof}
\noindent The growth assumptions imply $\eta\geq0$
and that $\mathcal J_j$ is finite on $\mathcal A_j$.

\vspace{0.2cm}

\noindent First,
$\mathcal A_j$ is nonempty, as
$(0,F^{j-1})\in\mathcal A_j$ and its infimum is finite as for every admissible pair $(v,F)$, the lower bound in
\eqref{eq:energy-growth} gives
\begin{equation}
\mathcal J_j(v,F)
\geq
\frac14\|v\|_{L^2}^2
+c_0\|F\|_{L^2}^2
-\frac12\|v^{j-1}\|_{L^2}^2
-c_1|\Om|.
\label{eq:coercivity}
\end{equation}

\vspace{0.2cm}

\noindent Let $(v_n,F_n)$ be a minimizing sequence.
By \eqref{eq:coercivity}, it is bounded in
$L^2(\Om;\mathbb R^d)\times
L^2(\Om;\mathbb R^{d\times d})$.
The constraint gives
\[
\|\nabla v_n\|_{L^2}
=
\frac1h\|F_n-F^{j-1}\|_{L^2}
\leq
\frac1h\left(\|F_n\|_{L^2}+\|F^{j-1}\|_{L^2}\right),
\]
and therefore $(v_n)$ is bounded in $W^{1,2}(\Om;\mathbb R^d)$.
Hence, after extracting a subsequence,
\[
\begin{aligned}
v_n&\rightharpoonup v^j
&&\text{in }W^{1,2}(\Om;\mathbb R^d),\\
F_n&\rightharpoonup F^j
&&\text{in }L^2(\Om;\mathbb R^{d\times d}).
\end{aligned}
\]
For every smooth periodic matrix test $\Phi$, the constraint
for the minimizing sequence reads
\[
\int_\Om(F_n-F^{j-1}):\Phi\,dx
=
-h\int_\Om v_n\cdot\operatorname{div}\Phi\,dx.
\]
Taking the limit $n\to\infty$ shows that the admissible class is sequentially weakly
closed in $W^{1,2}\times L^2$.

\vspace{0.2cm}

\noindent In particular,
\[
F^j-F^{j-1}=h\nabla v^j
\]
in the sense of distributions, and so
$(v^j,F^j)\in\mathcal A_j$.

\vspace{0.2cm}

\noindent For the fixed measurable function $\eta^{j-1}$,
the integrand in \eqref{eq:functional} is convex and continuous
in $(v,F)$. 
Consequently the functional is weakly lower
semicontinuous on the product $L^2$ space. Hence
\[
\mathcal J_j(v^j,F^j)
\leq\liminf_{n\to\infty}\mathcal J_j(v_n,F_n).
\]
and $(v^j,F^j)$ a minimizer.

\vspace{0.2cm}

\noindent To prove uniqueness, recall that the admissible
class is affine. The kinetic term is strictly convex
in $w$, while \eqref{eq:uniform-convexity} gives strict
convexity in $G$. More explicitly, for two admissible
pairs $(w_1,G_1)$ and $(w_2,G_2)$,
\begin{align*}
\mathcal J_j\left(\frac{w_1+w_2}{2},
                         \frac{G_1+G_2}{2}\right)
\leq
\frac12\mathcal J_j(w_1,G_1)
+\frac12\mathcal J_j(w_2,G_2)
-\frac18\|w_1-w_2\|_{L^2}^2
-\frac{\kappa}{8}\|G_1-G_2\|_{L^2}^2.
\end{align*}
Two minimizers must therefore coincide almost everywhere.
\end{proof}

\vspace{0.2cm}

\noindent\underline{\textbf{The discrete mechanical equations.}}
Let $\varphi\in W^{1,2}(\Om;\mathbb R^d)$.
The admissible variations are given by
\[
v_\varepsilon=v^j+\varepsilon\varphi,
\qquad
F_\varepsilon=F^j+\varepsilon h\nabla\varphi.
\]
Indeed, by \eqref{eq:kinematic-update},
\[
\begin{aligned}
F_\varepsilon-F^{j-1}
=F^j-F^{j-1}+\varepsilon h\nabla\varphi=h\nabla v^j+\varepsilon h\nabla\varphi
=h\nabla v_\varepsilon.
\end{aligned}
\]

\vspace{0.2cm}

\noindent The growth bound implies
\[
S(F^j,\eta^{j-1})
\in L^2(\Om;\mathbb R^{d\times d}).
\]
Indeed, for $|\varepsilon|\leq1$,
\[
\begin{aligned}
\left|
S(F^j+\varepsilon h\nabla\varphi,\eta^{j-1})
:h\nabla\varphi
\right|\leq
Ch\bigl(1+|F^j|+h|\nabla\varphi|+\eta^{j-1}\bigr)
|\nabla\varphi|
\in L^1(\Om).
\end{aligned}
\]
The derivative of the kinetic term is also dominated by
\[
\bigl(|v^j|+|v^{j-1}|+|\varphi|\bigr)|\varphi|
\in L^1(\Om).
\]
Moreover, both terms in the resulting
first variation define continuous linear functionals of
$\varphi\in W^{1,2}(\Om;\mathbb R^d)$.

\vspace{0.2cm}

\noindent Minimality of $(v^j,F^j)$ gives
\begin{equation}\label{EL}
\begin{aligned}
0
=
\left.\frac{d}{d\varepsilon}
\mathcal J_j(v_\varepsilon,F_\varepsilon)
\right|_{\varepsilon=0}=
\int_\Om
\left(
(v^j-v^{j-1})\cdot\varphi
+hS(F^j,\eta^{j-1}):\nabla\varphi
\right)\,dx.
\end{aligned}
\end{equation}
Hence the Euler--Lagrange equation leads to the discrete
momentum equation
\begin{equation}
\frac{v^j-v^{j-1}}h
=
\operatorname{div}S(F^j,\eta^{j-1})
\label{eq:momentum}
\end{equation}
in the sense of distributions.

\vspace{0.2cm}

\noindent More precisely, the equation holds in
$H^{-1}(\Om;\mathbb R^d)$, with
\[
\left\langle
\operatorname{div}S(F^j,\eta^{j-1}),\varphi
\right\rangle
=
-\int_\Om
S(F^j,\eta^{j-1}):\nabla\varphi\,dx.
\]
Since $(v^j-v^{j-1})/h\in L^2$ for fixed $h>0$,
the distributional divergence of the stress is represented
by an $L^2$ vector field.

\vspace{0.2cm}

\noindent For subsequent use, set
\[
S_*^j:=S(F^j,\eta^{j-1}).
\]
The discrete mechanical equations can then be written as
\begin{equation}\label{discretePDE}
\frac{F^j-F^{j-1}}h=\nabla v^j,
\qquad
\frac{v^j-v^{j-1}}h=\operatorname{div}S_*^j.
\end{equation}
If $F^0=\nabla y^0$ with
$y^0\in W^{1,2}(\Om;\mathbb R^d)$, define
\[
y^j=y^{j-1}+hv^j.
\]
Then \eqref{eq:kinematic-update} implies  that
$F^j=\nabla y^j$. Thus the construction preserves gradient
compatibility.

\vspace{0.2cm}

\noindent\underline{\textbf{The mechanical energy identity.}}
Next, we obtain the formula for the dissipation of the mechanical energy during 
the minimization step. Formally, this is achieved by multiplying \eqref{discretePDE}$_2$ by $v^j$
and then integrating by parts. To justify this  rigorously, observe that 
\[
e(F^{j-1},\eta^{j-1}),\ e(F^j,\eta^{j-1})
\in L^1(\Om).
\]
and that
\[
S(F^j,\eta^{j-1})\in L^2,
\qquad
F^j-F^{j-1}\in L^2,
\]
so their product belongs to $L^1$. Taking $\varphi=v^j\in W^{1,2}(\Om;\mathbb R^d)$ in \eqref{EL} and using
\eqref{eq:kinematic-update}, we infer that
\begin{equation}
\int_\Om
\left(
(v^j-v^{j-1})\cdot v^j
+S(F^j,\eta^{j-1}):(F^j-F^{j-1})
\right)\,dx=0.
\label{eq:mechanical-test}
\end{equation}

{
Next, note that the kinetic term satisfies the pointwise identity
\begin{equation}\label{eq:kinetic-identity}
(v^j-v^{j-1})\cdot v^j
=
\frac12|v^j|^2-\frac12|v^{j-1}|^2
+\frac12|v^j-v^{j-1}|^2.
\end{equation}
We employ the Bregman distance in the form
\begin{equation}\label{R}
\begin{aligned}
R^j
&:=  e \big ( F^{j-1} , \eta^{j-1} |  F^j , \eta^{j-1} \big )
\\
&= e(F^{j-1},\eta^{j-1})
-e(F^j,\eta^{j-1})
 - \frac{\del e}{\del F} (F^j,\eta^{j-1}):(F^{j-1} - F^j).
\end{aligned}
\end{equation}
Observe that this quantity constitutes the relative energy with the entropy variable fixed at the previous
step $\eta^{j-1}$. The uniform convexity  \eqref{eq:uniform-convexity} of $e(\cdot,\eta^{j-1})$ gives
\[
\begin{aligned}
R^j \ge \frac{\kappa}{2}|F^j-F^{j-1}|^2 
\end{aligned}
\]
Together \eqref{eq:kinetic-identity} and \eqref{R} yield the identity
\begin{equation}\label{relendiss}
\begin{aligned}
(v^j-v^{j-1}) &\cdot v^j
+S(F^j,\eta^{j-1}):(F^j-F^{j-1})
\\
&= \Big ( \frac12|v^j|^2 + e(F^j,\eta^{j-1}) \Big )  - \Big ( \frac12|v^{j-1}|^2 + e(F^{j-1} , \eta^{j-1}) \Big )
\\
&\qquad + \frac12|v^j-v^{j-1}|^2 + e \big ( F^{j-1} , \eta^{j-1} |  F^j , \eta^{j-1} \big ).
\end{aligned}
\end{equation}
Define the full mechanical dissipation density by
\begin{equation}\label{D}
D^j=\frac12|v^j-v^{j-1}|^2+R^j,
\end{equation}
which satisfies
\begin{equation}
D^j\geq
\frac12|v^j-v^{j-1}|^2
+\frac{\kappa}{2}|F^j-F^{j-1}|^2
\geq0.
\label{eq:dissipation-positive}
\end{equation}
Then  \eqref{eq:mechanical-test} yields the discrete mechanical energy balance
\begin{equation}\label{mechanical}
\begin{aligned}
\int_\Om
\left(\frac12|v^j|^2+e(F^j,\eta^{j-1})\right)\,dx
+\int_\Om D^j\,dx=
\int_\Om
\left(\frac12|v^{j-1}|^2
+e(F^{j-1},\eta^{j-1})\right)\,dx.
\end{aligned}
\end{equation}
In each mechanical step we lose exactly
$\int_\Om D^j\,dx$ from the total energy, which is converted
into internal energy in the next step.
}

\vspace{0.2cm}

\noindent\underline{\textbf{The entropy correction and energy balance.}}
For almost every $x\in\Om$, the map $s\mapsto e(F^j(x),s)$
is continuous, strictly increasing, and tends to infinity, since
$\partial_\eta e=\theta\geq\theta_0>0$.
In addition, $D^j+hr^j\geq0$, and so the equation
\begin{equation}\label{correction}
e(F^j,\eta^j)
=
e(F^j,\eta^{j-1})+D^j+hr^j
\end{equation}
determines a unique measurable $\eta^j\geq\eta^{j-1}\geq0$.

\vspace{0.2cm}

\noindent Moreover, $e(F^j,\eta^{j-1})$, $D^j$, and $r^j$
belong to $L^1$. The defining correction therefore shows
directly that $e(F^j,\eta^j)\in L^1$. The lower energy
bound then gives
\[
0\leq c_0|\eta^j|^2
\leq e(F^j,\eta^j)+c_1,
\]
and hence $\eta^j\in L^2(\Om)$.

\vspace{0.2cm}

\noindent Integrating \eqref{correction}, we obtain
\begin{equation}
\int_\Om e(F^j,\eta^j)\,dx
=
\int_\Om e(F^j,\eta^{j-1})\,dx
+\int_\Om D^j\,dx
+h\int_\Om r^j\,dx,
\label{eq:integrated-correction}
\end{equation}
and hence adding the kinetic energy and using \eqref{mechanical},
we infer that
\begin{equation}\label{eq:total-energy}
\begin{aligned}
\int_\Om
\left(\frac12|v^j|^2+e(F^j,\eta^j)\right)\,dx
=
\int_\Om
\left(\frac12|v^{j-1}|^2
+e(F^{j-1},\eta^{j-1})\right)\,dx+h\int_\Om r^j\,dx.
\end{aligned}
\end{equation}
We observe that there is no dissipation in the total energy
balance. The entropy correction restores precisely this amount
as internal energy. After the mechanical update and the entropy
correction, the total energy changes only by the prescribed heat
input. The dissipation still produces entropy, but it does not
cause a loss of total energy.

\vspace{0.2cm}

\noindent Summing \eqref{eq:total-energy} over $j$ gives
\begin{equation}\label{eq:summed-energy}
\begin{aligned}
\int_\Om
\left(\frac12|v^N|^2+e(F^N,\eta^N)\right)\,dx
=
\int_\Om
\left(\frac12|v^0|^2+e(F^0,\eta^0)\right)\,dx+
\int_0^{t_N}\int_\Om r(t,x)\,dx\,dt.
\end{aligned}
\end{equation}
In the notation of the preliminaries, the same calculation
up to any step $j$ reads
\begin{equation}\label{source-discrete-energy}
\int_\Om I(v^j,F^j,\eta^j)\,dx
=
E_0+\int_0^{t_j}\int_\Om r(t,x)\,dx\,dt
=
\mathcal E(t_j).
\end{equation}
Consequently, the lower bound in \eqref{eq:energy-growth}
implies
\begin{equation}\label{bounds}
\sup_{N\geq1}\ \max_{0\leq j\leq N}
\left(
\|v^j\|_{L^2}^2
+\|F^j\|_{L^2}^2
+\|\eta^j\|_{L^2}^2
\right)\leq C,
\end{equation}
where $C$ does not depend on $h$.

\vspace{0.2cm}

\noindent\underline{\textbf{The discrete entropy relation.}}
Having constructed the iterate $(v^j,F^j,\eta^j)$, we define
the averaged temperature
\begin{equation}
\theta_\#^j
:=
\int_0^1
\theta\left(
F^j,\eta^{j-1}+s(\eta^j-\eta^{j-1})
\right)\,ds.
\label{eq:secant-temperature}
\end{equation}
The fundamental theorem of calculus gives the identity
\begin{equation}
e(F^j,\eta^j)-e(F^j,\eta^{j-1})
=
\theta_\#^j(\eta^j-\eta^{j-1}),
\label{eq:secant-identity}
\end{equation}
while temperature growth gives
\[
\begin{aligned}
\theta_0\leq\theta_\#^j
\leq
C\int_0^1
\left(
1+|F^j|+\eta^{j-1}
+s(\eta^j-\eta^{j-1})
\right)\,ds\leq C(1+|F^j|+\eta^j).
\end{aligned}
\]

\vspace{0.2cm}

\noindent  The above bounds therefore imply
\[
\theta_\#^j\in L^2(\Om),
\qquad
\frac1{\theta_\#^j}\in L^\infty(\Om),
\qquad
\left\|\frac1{\theta_\#^j}\right\|_{L^\infty}
\leq\theta_0^{-1}.
\]
Since $\eta^j-\eta^{j-1}\in L^2$, the product in
\eqref{eq:secant-identity} belongs to $L^1$. By \eqref{correction} and
\eqref{eq:secant-identity}, we obtain
\begin{equation}\label{secant}
\theta_\#^j(\eta^j-\eta^{j-1})=D^j+hr^j.
\end{equation}
Dividing by $h\theta_\#^j$ and using the definition of
$D^j$, we conclude that
\begin{equation}\label{eq:entropy-production}
\begin{aligned}
\frac{\eta^j-\eta^{j-1}}h
&=
\frac{r^j}{\theta_\#^j}
+\frac{D^j}{h\theta_\#^j}
\\
&=
\frac{r^j}{\theta_\#^j}
+\frac{R^j}{h\theta_\#^j}
+\frac{|v^j-v^{j-1}|^2}{2h\theta_\#^j}.
\end{aligned}
\end{equation}

\vspace{0.2cm}

\noindent For each fixed $h>0$, all terms in
\eqref{eq:entropy-production} belong to $L^1$. Hence the identity and the resulting
inequality hold almost everywhere as well as in distributions.

\vspace{0.2cm}

\noindent The entropy increase consists of the contribution
from the prescribed heat supply and from converting the
mechanical dissipation into internal energy. Thus the correction
restores the energy lost during the mechanical update while
increasing the entropy. The nonnegativity of the terms on the
right-hand side gives the discrete entropy production inequality
\begin{equation}
\frac{\eta^j-\eta^{j-1}}h
\geq\frac{r^j}{\theta_\#^j}\geq0.
\label{eq:entropy-inequality}
\end{equation}
In particular, the last inequality may be multiplied
by any $\zeta\in L^\infty(\Om)$, $\zeta \ge 0$ to yield
\[
\int_\Om
\frac{\eta^j-\eta^{j-1}}h\,\zeta\,dx
\geq
\int_\Om
\frac{r^j}{\theta_\#^j}\,\zeta\,dx
\geq0.
\]

\vspace{0.2cm}

{
\noindent
{\bf Summary:} It is instructive to summarize the equations that the iterates $(v^j , F^j, \eta^j)$ satisfy: By construction
they satisfy the relation \eqref{discretePDE} and \eqref{correction}:
\begin{equation}\label{iterates}
\begin{aligned}
\frac{F^j-F^{j-1}}h &=\nabla v^j,
\\
\frac{v^j-v^{j-1}}h &=\operatorname{div}S (F^j, \eta^{j-1} )
\\[5pt]
e(F^j,\eta^j) -  e(F^j,\eta^{j-1}) &= D^j+hr^j
\end{aligned}
\end{equation}
From the first two equations in \eqref{iterates} one derives the equation
\begin{equation}\label{mechdiss}
\Big ( \frac12|v^j|^2 + e(F^j,\eta^{j-1}) \Big )  - \Big ( \frac12|v^{j-1}|^2 + e(F^{j-1}\eta^{j-1}) \Big ) + D^j 
= h \diver \big( v^j \cdot  S (F^j, \eta^{j-1}  )  \big)
\end{equation}
Formally, this is easily done  using  the identity \eqref{relendiss}; to do it rigorously one has to validate the product rule
$$
\diver \big( v^j \cdot  S (F^j, \eta^{j-1}  ) \big)  = S (F^j, \eta^{j-1} ) : \nabla v^j  + v^j \cdot \diver \big( S (F^j, \eta^{j-1}  ) \big ) \, .
$$
The latter can be justified within the available regularity for the iterates, if we take into account that the second equation in \eqref{iterates}
implies that $\diver S(F^j,\eta^{j-1}) \in L^2$. 
Using \eqref{D} and \eqref{eq:secant-identity},
the third equation in \eqref{iterates} and \eqref{mechdiss} yield :
\begin{align}
\label{energycons}
\frac{\Big ( \frac12|v^j|^2 + e(F^j,\eta^j ) \Big )  - \Big ( \frac12|v^{j-1}|^2 + e(F^{j-1} , \eta^{j-1}) \Big )}{h}
&= \diver \big( v^j \cdot  S (F^j, \eta^{j-1}  ) \big ) + r^j \, , 
\\[5pt]
\label{entropyprod}
\frac{\eta^j-\eta^{j-1}}h
&= \frac{r^j}{\theta_\#^j} +\frac{D^j}{h\theta_\#^j} \, , 
\end{align}
the discrete variants of the conservation of energy and entropy production identities, respectively.
}

\section{Construction of a dissipative measure-valued solution}\label{sec:construction}

Let $\mathcal M_+(\Om)$ denote the finite nonnegative
Radon measures on $\Om$.

\begin{theorem}\label{thm}
Assume \emph{(H1)--(H6)}, and prescribe
$r\in L^1(Q_T)$, $r\geq0$.
For every $T>0$ and initial data
\[
 y^0\in W^{1,2}(\Om;\R^d),\quad
 F^0=\nabla y^0,\quad
 v^0\in L^2(\Om;\R^d),\quad
 \eta^0\in L^2(\Om),\quad\eta^0\geq0,
\]
there exists a dissipative measure-valued solution in the
sense of Definition \ref{defDMV}.

The discrete construction is well defined at every step,
its mechanical minimizer is unique, and
\begin{equation}\label{energy-discrete}
 \int_\Om I(v^j,F^j,\eta^j)\dd x=\mathcal E(t_j),
 \qquad
 \eta^j\geq\eta^{j-1}\geq0.
\end{equation}
For every sequence $h=T/N\to 0$, there is a
subsequence generating a
Young measure $\nu_{t,x}$ on $\mathcal Z$, a finite
nonnegative Radon measure $\sigma$ on
$Q_T=(0,T)\times\Om$, and a weak-* measurable family
$\gamma_t\in\mathcal M_+(\Om)$  with the following properties.

The second moment of the Young measure satisfies
\begin{equation}\label{moments}
 \operatorname*{ess\,sup}_{t\in(0,T)}
 \int_\Om\nuav{|F|^2+|v|^2+\eta^2}\dd x<\infty.
\end{equation}
Set
\[
 \bar F=\nuav{F},\qquad
 \bar v=\nuav{v},\qquad
 \bar\eta=\nuav{\eta},\qquad
 \bar S=\nuav{S(F,\eta)}.
\]
For every smooth periodic test function, compactly supported
in time in $[0,T)$, one has
\begin{align}
 \int_{\QT}
 \big(\bar F:\partial_t\Phi
      -\bar v\cdot\diver\Phi\big)\dd x\dd t
 +\int_\Om F^0:\Phi(0)\dd x&=0,
 \label{weakF}\\
 \int_{\QT}
{ \big(\bar v\cdot\partial_t\varphi
      - \bar S:\nabla\varphi\big)\dd x\dd t }
 +\int_\Om v^0\cdot\varphi(0)\dd x&=0,
 \label{weakv}\\
 \int_{\QT}\bar\eta\,\partial_t\zeta\dd x\dd t+
 \int_\Om\eta^0\zeta(0)\dd x +\int_{Q_T}r\nuav{1/\theta(F,\eta)}\zeta\dd x\dd t
 +\int_{Q_T}\zeta\dd\sigma&=0.
 \label{weaketa}
\end{align}
Here $\Phi$, $\varphi$ and $\zeta$ are matrix-,
vector- and scalar-valued, respectively; $\zeta$ need
not be nonnegative. Thus
\[
 \partial_t\bar F=\nabla\bar v,\qquad
 \partial_t\bar v=\diver\bar S,\qquad
 \partial_t\bar\eta
 =r\nuav{1/\theta(F,\eta)}+\sigma
\]
in the sense of distributions.

Up to redefining the barycentres on a negligible set in time, the mechanical barycentres are weakly continuous in $L^2$ with $\bar F(0)=F^0$,
$\bar v(0)=v^0$. The entropy is nondecreasing and
\begin{equation}\label{entropytrace}
 \lim_{t\to 0}
 \|\bar\eta(t)-\eta^0\|_{L^1(\Om)}=0.
\end{equation}
For almost every $t\in(0,T)$,
\begin{equation}\label{energy-limit}
 \int_\Om\nuav{I(v,F,\eta)}\dd x+\gamma_t(\Om)
 =\mathcal E(t),
\end{equation}
while
\begin{equation}\label{lag}
 \|\eta_h-\eta_h^-\|_{L^2(\QT)}\leq Ch^{1/2},
 \qquad
 S_h^*\rightharpoonup\bar S
 \quad\hbox{in }L^2(\QT).
\end{equation}
Each $F^j$ is a gradient, and the motion satisfies
\[
 y\in W^{1,\infty}(0,T;L^2)
       \cap L^\infty(0,T;W^{1,2}).
\]
The concentration measure is
$\gamma(\dd t,\dd x)=\gamma_t(\dd x)\dd t$,
so \eqref{energy-limit} implies \eqref{DMVenergy}.
\end{theorem}

\begin{proof}

\noindent\underline{Step 1: A priori estimates.} Take $h=T/N$.
On $(t_{j-1},t_j]$, put
\begin{equation}\label{interpolants}
 \begin{gathered}
 u_h=(v_h,f_h,\eta_h)=(v^j,F^j,\eta^j),\\
 \eta_h^-=\eta^{j-1},\quad\text{and}\quad
 S_h^*=S(f_h,\eta_h^-).
 \end{gathered}
\end{equation}
For any discrete quantity $a^j$, write $A_h$ for
its affine interpolant. In particular,
$a^j=v^j,F^j,\eta^j$ give $V_h,F_h,H_h$.
Thus
\[
 A_h(t)
 =a^{j-1}+\frac{t-t_{j-1}}h(a^j-a^{j-1}),
 \qquad t\in[t_{j-1},t_j].
\]

\noindent In addition, we define the triples
\[
 U_h=(V_h,F_h,H_h),
 \qquad
 W_h=(v_h,F_h,H_h).
\]
Note that both $U_h$ and $W_h$ will generate the same Young measure as the
piecewise constant interpolants. The distinction is that the velocity
in $W_h$ is exactly the time derivative of the affine
deformation. Throughout the proof, $C$ denotes a constant independent of $h$ and we use the notation
$L^\infty L^2=L^\infty(0,T;L^2(\Om))$ and
$H^{-1}=(W^{1,2}(\Om))^*$.

\vspace{0.2cm}

\noindent We first estimate the entropy increments and use
these estimates, together with uniform convexity in the
deformation gradient, to show that the piecewise constant
and affine interpolants differ by quantities tending to
zero strongly in $L^2(Q_T)$. Write $\Delta\eta^j=\eta^j-\eta^{j-1}\geq0$, then
pointwise monotonicity gives
\[
 \sum_{j=1}^{N}|\Delta\eta^j|^2
 \leq\left(\sum_{j=1}^{N}\Delta\eta^j\right)^2
 =|\eta^{N}-\eta^0|^2
 \leq|\eta^{N}|^2.
\]
Consequently,
\begin{equation}\label{lag-estimate}
\begin{aligned}
 \|\eta_h-\eta_h^-\|_{L^2(Q_T)}^2
 &=h\sum_{j=1}^{N}\|\Delta\eta^j\|_2^2
 \leq Ch,
\\
 \|H_h-\eta_h\|_{L^2(Q_T)}^2
 &=\frac h3\sum_{j=1}^N\|\Delta\eta^j\|_2^2
 \leq Ch.
\end{aligned}
\end{equation}
and
\begin{equation}\label{lagL1}
 \|\eta_h-\eta_h^-\|_{L^1(Q_T)}
 =
 h\int_\Om(\eta^{N}-\eta^0)\dd x
 \leq Ch.
\end{equation}

\vspace{0.2cm}

\noindent Define $D_h=D^j$ on each cell, then by (H6) and
\eqref{secant} we infer that
\[
 0\leq D_h
 \leq C(1+|f_h|+\eta_h)(\eta_h-\eta_h^-).
\]
The first factor is bounded in $L^2(Q_T)$ by
\eqref{bounds}, while Cauchy--Schwarz and
\eqref{lag-estimate} imply
\begin{equation}\label{Dsmall}
\begin{aligned}
 \int_{Q_T}D_h\dd x\dd t
 &= h \sum_{j=1}^N \int_\Om D^j \dd x
 \leq Ch^{1/2}.
 \end{aligned}
\end{equation}
Uniform convexity yields
\[
 D^j\geq
 \frac{\kappa}{2}|F^j-F^{j-1}|^2
 +\frac12|v^j-v^{j-1}|^2.
\]
By a direct calculation, using the definitions of the piecewise contant $f_h, v_h $ and picewise linear $F_h, V_h$
interpolations, we obtain 
\begin{equation}\label{strongconsistency}
 \|F_h-f_h\|_{L^2(Q_T)}^2
 +\|V_h-v_h\|_{L^2(Q_T)}^2
 \leq Ch^{1/2}.
\end{equation}
In particular, these differences converge strongly in
$L^2$. 


\vspace{0.4cm}

\noindent\underline{Step 2: Young-measure generation.} We are now ready to pass to the limit in the discrete equations. The estimates above give
\[
\sup_h\left(
\|u_h\|_{L^\infty L^2}
+\|U_h\|_{L^\infty L^2}
+\|W_h\|_{L^\infty L^2}
+\|\eta_h^-\|_{L^\infty L^2}
\right)\leq C.
\]
Moreover, \eqref{lag-estimate} and
\eqref{strongconsistency} imply
\[
\|U_h-u_h\|_{L^2(Q_T)}
+\|W_h-u_h\|_{L^2(Q_T)}
\longrightarrow0.
\]
The energy bound above provides
weak compactness of the states, while the entropy-increment
estimate and uniform convexity identify the Young measures generated
by the different interpolants.

\vspace{0.2cm}

\noindent By \eqref{bounds}
there exists a subsequence, not relabelled, generating a Young measure $\nu_{t,x}$ on $\mathcal Z$
satisfying \eqref{moments}. Writing
\[
(\bar v,\bar F,\bar\eta)
=
\nuav{(v,F,\eta)},
\]
we have
\[
u_h\overset{*}{\rightharpoonup}
(\bar v,\bar F,\bar\eta)
\qquad\text{in }L^\infty(0,T;L^2(\Om)).
\]
Moreover, for every continuous function
$g:\mathcal Z\to\R$ with at most linear growth,
\[
g(u_h)\rightharpoonup\nuav{g}
\qquad\text{in }L^1(Q_T).
\]

\noindent By \eqref{lag-estimate} and
\eqref{strongconsistency},
\[
\|U_h-u_h\|_{L^2(Q_T)}
+
\|W_h-u_h\|_{L^2(Q_T)}
\longrightarrow0.
\]
Hence $U_h$ and
$W_h$ generate the same Young measure $\nu$.
Their uniform $L^\infty L^2$ bounds also give the same
weak-* limit $(\bar v,\bar F,\bar\eta)$.

\vspace{0.2cm}

\noindent We verify that $\nu$ is generated by a sequence
with the compatibility and regularity required in
Definition \ref{defDMV}. Recall that
\[
y^j=y^{j-1}+hv^j,
\qquad
\nabla y^j=F^j,
\]
and define
\[
Y_h(t)=y^0+\int_0^t v_h(s)\dd s.
\]
On each time cell, $Y_h$ is the affine interpolant of
$y^{j-1}$ and $y^j$, and consequently
\[
\partial_tY_h=v_h,
\qquad
\nabla Y_h=F_h.
\]
In particular,
\[
\sup_{t\in[0,T]}\|Y_h(t)\|_{L^2}
\leq \|y^0\|_{L^2}
+T\|v_h\|_{L^\infty L^2}\leq C.
\]
Combining this with the bound on $F_h$ gives
\begin{equation}\label{generator-bounds}
\begin{gathered}
\|Y_h\|_{L^\infty(0,T;W^{1,2})}
+\|\partial_tY_h\|_{L^\infty L^2}
+\|H_h\|_{L^\infty L^2}\leq C,\\
\|\partial_t\nabla Y_h\|_{L^\infty(0,T;H^{-1})}
\leq C.
\end{gathered}
\end{equation}
For the last estimate, we use
$\partial_t\nabla Y_h=\nabla v_h$ in the sense of distributions and
\[
|\langle\nabla v_h(t),\Phi\rangle|
=
\left|\int_\Om v_h(t)\cdot\diver\Phi\dd x\right|
\leq C\|v_h(t)\|_{L^2}\|\Phi\|_{W^{1,2}},
\]
for all $\Phi\in W^{1,2}(\Om;\R^{d\times d})$.  Since $W_h=(\partial_tY_h,\nabla Y_h,H_h)$
generates $\nu$ and $H_h\geq0$, estimate
\eqref{generator-bounds} verifies the requirements for the generating sequence
in Definition \ref{defDMV}.

\vspace{0.2cm}

\noindent The corresponding energy measures are consistent
as well. By (H5)--(H6), integration along the segment
joining $z,z'\in\mathcal Z$ gives
\[
|I(z)-I(z')|
\leq C(1+|z|+|z'|)|z-z'|.
\]
Hence Cauchy--Schwarz and the preceding bounds imply
\begin{equation}\label{energyconsistency}
\|I(W_h)-I(u_h)\|_{L^1(Q_T)}
+
\|I(U_h)-I(u_h)\|_{L^1(Q_T)}
\longrightarrow0.
\end{equation}
This estimate will allow us to use the exact energy
balance of $u_h$ when identifying the energy measure
associated with the admissible generating sequence.

\vspace{0.2cm}

\noindent
{The equations  \eqref{iterates} and \eqref{energycons} for the iterates imply that
the approximate solutions satisfy the conservation equations:
\begin{align}
\del_t F_h  &= \nabla v_h
\\
\del_t V_h  &= \diver S(f_h , \eta_h^-)
\\
\del_t {\mathcal J}_h &= \diver \big (  v_h \cdot S(f_h , \eta_h^- ) \big ) + r_h .
\end{align}
Above ${\mathcal J}_h$ stands for the piecewise linear interpolation,
$$
{\mathcal J}_h = I(U^{j-1}) + \frac{t - t_{j-1}}{h} \big(I(U^j) - I(U^{j-1}) \big) \quad \mbox{on the cell $(t_{j-1}, t_j]$ where $U^j = (v^j, F^j, \eta^j)$} \, ,
$$
and we recall that $S(f_h , \eta_h^-) = S(F^j , \eta^{j-1})$ and $r_h = r^j$ on each cell $(t_{j-1}, t_j]$.
The equation for the entropy production will be undertaken later.
}

\vspace{0.2cm}
\noindent
We prove an integrated form of the energy identity.

\begin{lemma}
For every
$\varphi\in C_c^1([0,T))$,
\begin{equation}\label{affineenergybalance}
E_0\varphi(0)
+\int_{Q_T}\varphi'(t)I(U_h)\dd x\dd t
+\int_{Q_T}r\varphi\dd x\dd t=o(1).
\end{equation}
\end{lemma} 

\begin{proof}
{Set first
\[
\begin{aligned}
R(t)=\int_\Om r(t,x)\dd x \, ,
\\
\mathcal E (t) = E_0 + \int_0^t R (s) \dd s \, , 
\end{aligned}
\]
and note that, since $\mathcal E\in W^{1,1}(0,T)$,
$\mathcal E'=R$ a.e., and
$\mathcal E(0)=E_0$, integration by parts gives the identity
\[
E_0\varphi(0)
+\int_0^T\varphi'(t)\mathcal E(t)\dd t
+\int_{Q_T}r\varphi\dd x\dd t=0 \, .
\]
}

\vspace{0.2cm}
\noindent
On the other hand,  \eqref{source-discrete-energy} implies
\[
\int_\Om I(u_h(t,x))\dd x=\mathcal E(t_j)
\qquad\text{for }t\in(t_{j-1},t_j].
\]
Thus the absolute value of the left-hand side of
\eqref{affineenergybalance} is at most
\[
\|\varphi'\|_\infty
\left(
\|I(U_h)-I(u_h)\|_{L^1(Q_T)}
+
\sum_{j=1}^N\int_{t_{j-1}}^{t_j}
|\mathcal E(t_j)-\mathcal E(t)|\dd t
\right).
\]
The first term tends to zero by
\eqref{energyconsistency}, while for the second $R\geq0$ gives
\[
\begin{aligned}
\sum_{j=1}^N\int_{t_{j-1}}^{t_j}
|\mathcal E(t_j)-\mathcal E(t)|\dd t
&=
\sum_{j=1}^N\int_{t_{j-1}}^{t_j}
\int_t^{t_j}R(s)\dd s\dd t\\
&=
\sum_{j=1}^N\int_{t_{j-1}}^{t_j}
(s-t_{j-1})R(s)\dd s\leq h\|r\|_{L^1(Q_T)}.
\end{aligned}
\]
This proves \eqref{affineenergybalance}.
\end{proof}

\vspace{0.2cm}

\noindent Since $\eta_h\geq\eta^0$ almost everywhere,
we claim that
\begin{equation}\label{support}
\nu_{t,x}\{(v,F,\eta)\in\mathcal Z:
\eta\geq\eta^0(x)\}=1,
\quad\text{for almost every }(t,x)\in Q_T.
\end{equation}
\noindent Indeed, consider the sequence
$(u_h,\eta^0)$. Since the second component is independent
of $h$, the standard properties of Young measures imply
that this sequence generates
$\nu_{t,x}\otimes\delta_{\eta^0(x)}$.
Apply the representation theorem to the bounded continuous
function
\[
G(v,F,\eta,a)=\min\{1,(a-\eta)_+\},
\qquad s_+=\max\{s,0\}.
\]
As $\eta_h\geq\eta^0$ almost everywhere, we have
$G(u_h,\eta^0)=0$, and hence
\[
0
=
\lim_{h\to0}\int_{Q_T}G(u_h,\eta^0)\dd x\dd t
=
\int_{Q_T}
\left\langle\nu_{t,x},
\min\{1,(\eta^0(x)-\eta)_+\}
\right\rangle\dd x\dd t.
\]
The integrand is nonnegative and therefore vanishes
almost everywhere. Since
$\min\{1,(\eta^0(x)-\eta)_+\}>0$ precisely when
$\eta<\eta^0(x)$, it follows that
\[
\nu_{t,x}\{\eta<\eta^0(x)\}=0
\qquad\text{for almost every }(t,x)\in Q_T.
\]
This proves \eqref{support} and, in particular,
$\bar\eta=\nuav{\eta}\geq\eta^0$ almost everywhere.

\vspace{0.4cm}

\noindent\underline{Step 3: The mechanical equations.} We next identify the weak limit of the stress $S_h^*=S(f_h,\eta_h^-)$. By (H5) and the uniform energy bound,
\[
\sup_h\left(
\|S(f_h,\eta_h^-)\|_{L^\infty L^2}
+
\|S(f_h,\eta_h)\|_{L^\infty L^2}
\right)\leq C.
\]
We claim that
\[
\|S(f_h,\eta_h^-)-S(f_h,\eta_h)\|_{L^1(Q_T)}
\longrightarrow0.
\]
For $M>0$, set
\begin{equation}\label{definitionA}
A_{h,M}
=
\{(t,x)\in Q_T:|f_h|+\eta_h+\eta_h^->M\}.
\end{equation}
The second-moment estimates and Chebyshev's inequality imply
$|A_{h,M}|\leq C/M^2$, and therefore
\[
\begin{aligned}
&\int_{A_{h,M}}
\bigl(|S(f_h,\eta_h^-)|+|S(f_h,\eta_h)|\bigr)
\dd x\dd t
\\
&\qquad \leq
|A_{h,M}|^{1/2}
\left(
\|S(f_h,\eta_h^-)\|_{L^2(Q_T)}
+\|S(f_h,\eta_h)\|_{L^2(Q_T)}
\right)
\leq\frac{C}{M}.
\end{aligned}
\]
On $Q_T\setminus A_{h,M}$, the arguments lie in a fixed
compact set. Uniform continuity of $S$ and
\eqref{lag-estimate} imply convergence of their difference
to zero in measure. Since the stresses are bounded
on that compact set, the difference converges to zero
in $L^1$ on the complement of $A_{h,M}$.
Taking first $h\to0$ and then $M\to\infty$ proves the claim.

\vspace{0.2cm}

\noindent The linear growth of $S$ implies that
\[
S(f_h,\eta_h)\rightharpoonup
\bar S:=\nuav{S(F,\eta)}
\qquad\text{in }L^1(Q_T),
\]
which in turn implies
$S_h^*=S(f_h,\eta_h^-)\rightharpoonup\bar S$
in $L^1(Q_T)$, by the previous estimates.
Since $(S_h^*)$ is bounded in
$L^\infty(0,T;L^2(\Om;\R^{d\times d}))$,
we may extract a further subsequence such that
\[
S_h^*\overset{*}{\rightharpoonup}\bar S
\qquad\text{in }
L^\infty(0,T;L^2(\Om;\R^{d\times d})),
\]
and, in particular,
\[
S_h^*\rightharpoonup\bar S
\qquad\text{in }L^2(Q_T;\R^{d\times d}),
\]
which proves the stress convergence in \eqref{lag}.


\vspace{0.2cm}

\noindent We next pass to the limit in the discrete
kinematic and momentum equations to establish
\eqref{weakF} and \eqref{weakv}.
To this end, recall that the affine interpolants satisfy
\[
\partial_tF_h=\nabla v_h,
\qquad
\partial_tV_h=\diver S_h^*,
\qquad
F_h(0)=F^0,
\qquad
V_h(0)=v^0,
\]
where these identities hold in $L^\infty(0,T;H^{-1})$.
In particular,
\[
\|\partial_tF_h\|_{L^\infty H^{-1}}
+
\|\partial_tV_h\|_{L^\infty H^{-1}}
\leq C.
\]
On the $j$th cell,
\[
F_h(t)-f_h(t)=-(t_j-t)\nabla v^j,
\qquad
V_h(t)-v_h(t)=-(t_j-t)\diver S_*^j,
\]
as equalities in $H^{-1}$, and so
\[
\|F_h-f_h\|_{L^\infty H^{-1}}
+
\|V_h-v_h\|_{L^\infty H^{-1}}
\leq Ch.
\]

\vspace{0.2cm}

\noindent Let
\[
\Phi\in C_c^1([0,T)\times\Om;\R^{d\times d}),
\qquad
\varphi\in C_c^1([0,T)\times\Om;\R^d),
\]
be periodic test functions in space which vanish near $T$.
Integration by parts gives
\[
\int_{Q_T}
\bigl(F_h:\partial_t\Phi-v_h\cdot\diver\Phi\bigr)
\dd x\dd t
+\int_\Om F^0:\Phi(0)\dd x=0
\]
and
\[
\int_{Q_T}
{  \bigl(V_h\cdot\partial_t\varphi - S_h^*:\nabla\varphi\bigr) }
\dd x\dd t
+\int_\Om v^0\cdot\varphi(0)\dd x=0.
\]
All terms are integrable by the $L^\infty L^2$ bounds.
Passing to the limit proves \eqref{weakF} and
\eqref{weakv}, with the prescribed initial terms.

\vspace{0.2cm}

\noindent Consequently,
\[
\partial_t\bar F=\nabla\bar v,
\qquad
\partial_t\bar v=\diver\bar S,
\]
as distributions and both time derivatives belong to
$L^\infty(0,T;H^{-1})$.
Since also $\bar F,\bar v\in L^\infty L^2$, this gives
\[
\bar F\in C_{\mathrm w}([0,T];L^2(\Om;\R^{d\times d})),
\qquad
\bar v\in C_{\mathrm w}([0,T];L^2(\Om;\R^d)),
\]
i.e. $\bar F$ and $\bar v$ are weakly continuous in time with values in $L^2$.
\vspace{0.2cm}

\noindent It remains to identify a motion $y$ such that
$\partial_t y=\bar v$ and $\nabla y=\bar F$.
To this end, define
\[
y(t)=y^0+\int_0^t\bar v(s)\dd s,
\quad\text{in }L^2(\Om;\R^d).
\]
Then $y\in W^{1,\infty}(0,T;L^2)$,
$\partial_t y=\bar v$, and $y(0)=y^0$ in $L^2$.
The kinematic equation and its initial condition imply
\[
\nabla y(t)
=
F^0+\int_0^t\nabla\bar v(s)\dd s
=
\bar F(t),
\quad\text{in }H^{-1},
\]
and since $\bar F\in L^\infty L^2$, it follows that
\[
y\in W^{1,\infty}(0,T;L^2)
\cap L^\infty(0,T;W^{1,2}),
\quad
\nabla y=\bar F.
\]

\vspace{0.4cm}

\noindent\underline{Step 4: Entropy production.} We next establish the entropy relation
\eqref{weaketa}. Set $r_h=r^j$ and
$\theta_{\#,h}=\theta_\#^j$ on each cell
$(t_{j-1},t_j]$. 
Approximating $r$ in $L^1(Q_T)$ by continuous functions
and observing that their time averages converge uniformly,
we obtain
\[
r_h\longrightarrow r
\qquad\text{strongly in }L^1(Q_T).
\]
We next compare the averaged temperature with the temperature
evaluated at the corrected state. By definition,
\[
\theta_{\#,h}
=
\int_0^1
\theta\bigl(f_h,\eta_h^-+s(\eta_h-\eta_h^-)\bigr)\dd s.
\]
Outside the set $A_{h,M}$ introduced in \eqref{definitionA}, both the
endpoint $(f_h,\eta_h)$ and the entire entropy segment
belong to
\[
K_M=\{(F,a):|F|\leq M,\ 0\leq a\leq M\}.
\]
Here we used $0\leq\eta_h^-\leq\eta_h$.

\vspace{0.2cm}

\noindent Let $\omega_M$ be a nondecreasing modulus of
continuity of $\theta$ on $K_M$. Since
\[
|\eta_h^-+s(\eta_h-\eta_h^-)-\eta_h|
=(1-s)|\eta_h-\eta_h^-|,
\]
we have
\[
|\theta_{\#,h}-\theta(f_h,\eta_h)|
\leq\omega_M(|\eta_h-\eta_h^-|)
\qquad\text{on }Q_T\setminus A_{h,M}.
\]
For fixed $\varepsilon>0$ and $M>0$, choose $\rho>0$
such that $\omega_M(\rho)<\varepsilon$.
Then \eqref{lag-estimate} gives
\[
\begin{aligned}
\bigl|\{|\theta_{\#,h}-\theta(f_h,\eta_h)|>\varepsilon\}\bigr|
\leq
|A_{h,M}|+|\{|\eta_h-\eta_h^-|>\rho\}|&\leq
\frac{C}{M^2}
+\frac{\|\eta_h-\eta_h^-\|_{L^2(Q_T)}^2}{\rho^2}\\
&\leq\frac{C}{M^2}+\frac{Ch}{\rho^2}.
\end{aligned}
\]
Letting first $h\to0$ and then $M\to\infty$ proves
\[
\theta_{\#,h}-\theta(f_h,\eta_h)\longrightarrow0
\qquad\text{in measure on }Q_T.
\]

\vspace{0.2cm}

\noindent By (H4), both temperatures are at least
$\theta_0$. Hence
\[
\left|
\frac1{\theta_{\#,h}}
-\frac1{\theta(f_h,\eta_h)}
\right|
\leq
\frac1{\theta_0^2}
|\theta_{\#,h}-\theta(f_h,\eta_h)|,
\]
and so the above difference converges
to zero in measure and is uniformly bounded.
As $(v,F,\eta)\mapsto1/\theta(F,\eta)$ is bounded and
continuous, it follows that
\[
\frac1{\theta_{\#,h}}
\overset{*}{\rightharpoonup}
\nuav{1/\theta(F,\eta)}
\qquad\text{in }L^\infty(Q_T).
\]
For every fixed $\zeta\in L^\infty(Q_T)$, we may write
\[
\begin{aligned}
&\int_{Q_T}\zeta
\left(
\frac{r_h}{\theta_{\#,h}}
-r\nuav{1/\theta(F,\eta)}
\right)\dd x\dd t
\\
&\qquad =
\int_{Q_T}
\frac{\zeta(r_h-r)}{\theta_{\#,h}}\dd x\dd t
+
\int_{Q_T}\zeta r
\left(
\frac1{\theta_{\#,h}}
-\nuav{1/\theta(F,\eta)}
\right)\dd x\dd t.
\end{aligned}
\]
The first integral on the right-hand side is bounded by
$\theta_0^{-1}\|\zeta\|_\infty\|r_h-r\|_{L^1}$, while the second tends to zero by weak-* convergence, since
$\zeta r\in L^1(Q_T)$. Therefore
\begin{equation}\label{eq:source}
\frac{r_h}{\theta_{\#,h}}
\rightharpoonup
r\nuav{1/\theta(F,\eta)}
\qquad\text{in }L^1(Q_T).
\end{equation}

\vspace{0.2cm}

\noindent The remaining part of the discrete entropy
production is represented by the nonnegative measures
\[
\dd\sigma_h
=
\frac{D_h}{h\theta_{\#,h}}\dd x\dd t
\]
on $(0,T)\times\Om$. The discrete entropy identity \eqref{eq:entropy-production} gives
\[
\partial_tH_h
=
\frac{r_h}{\theta_{\#,h}}
+\frac{D_h}{h\theta_{\#,h}}
\qquad\text{almost everywhere in }Q_T.
\]
Integrating the above identity we find
\[
\begin{aligned}
\sigma_h(Q_T)
\leq
\int_{Q_T}\partial_tH_h\dd x\dd t=
\int_\Om(\eta^N-\eta^0)\dd x
\leq C.
\end{aligned}
\]
Thus, up to a further subsequence, $\sigma_h$ converges weakly-$\star$
in the sense measures on $Q_T$ to a limiting measure $\sigma$ .

\vspace{0.2cm}

\noindent Let $\zeta\in C_c^1((0,T)\times\Om)$. Integration by parts then gives
\[
\begin{aligned}
-\int_{Q_T}H_h\partial_t\zeta\dd x\dd t
&=
\int_{Q_T}\zeta\frac{r_h}{\theta_{\#,h}}\dd x\dd t
+\int_{Q_T}\zeta\dd\sigma_h.
\end{aligned}
\]
By the weak convergence of $H_h$, \eqref{eq:source}, and the weak-$\star$ convergence of $\sigma_h$, we infer that
\[
\begin{aligned}
-\int_{Q_T}\bar\eta\,\partial_t\zeta\dd x\dd t
=
\int_{Q_T}
r\nuav{1/\theta(F,\eta)}\zeta\dd x\dd t+
\int_{Q_T}\zeta\dd\sigma.
\end{aligned}
\]
To prove \eqref{weaketa} we need to test against $\zeta\in C_c^1([0,T)\times\Om)$ which we postpone until Step 6. Then, since $\sigma\geq0$, it also
gives the entropy inequality in Definition \ref{defDMV}.

\vspace{0.2cm}

\noindent Also, note that in $\QT$,
\[
\partial_t\bar\eta
=
r\nuav{1/\theta(F,\eta)}+\sigma\geq0,
\]
in the sense of distributions on $Q_T$. This implies that $\bar\eta$ admits a representative
nondecreasing in time for almost every $x$ and the support property \eqref{support} also implies
$\bar\eta(t)\geq\eta^0$.

\vspace{0.4cm}

{\color{black}
\noindent\underline{Step 5: Energy balance.}
Set $\widetilde I=I+C_*\geq0$ and note $\widetilde I \ge 0$ satisfies \eqref{shiftgrowth}. The uniform energy bound gives
\[
\sup_h\operatorname*{ess\,sup}_{t\in(0,T)}
\int_\Omega\widetilde I(u_h(t,x))\,dx<\infty.
\]
After extracting a further subsequence, there is a weak-* measurable
family of nonnegative measures $\mu_t$ such that
\[
\widetilde I(u_h(t,\cdot))\,dx
\overset{*}{\rightharpoonup}\mu_t
\quad\hbox{in }L^\infty_{w*}(0,T;\mathcal M(\Omega)).
\]
Here convergence is tested against $L^1(0,T;C(\Omega))$.

\vspace{0.2cm}
\noindent
Young measures are not directly defined for $\tilde I$ because of its quadratic growth. 
Using the positivity of $I(U)$, following \cite[App A.2]{DST2012}, consider for $M> 0$
the bounded continuous truncation
$\widetilde I_M=\min\{\widetilde I,M\}$ whose weak limits are represented by the Young measure.
Define 
$$\langle \nu, \widetilde I \rangle := \lim_{M \to \infty} \langle \nu, \widetilde I_M \rangle $$
using the
monotone convergence theorem. Testing with nonnegative functions and then letting $M\to\infty$
shows that
\[
\mu_t\geq\langle\nu_{t,x},\widetilde I\rangle\,dx
\quad\hbox{for almost every }t.
\]
Define
\[
\gamma_t=\mu_t-\langle\nu_{t,x},\widetilde I\rangle\,dx\geq0,
\qquad \gamma(dt,dx)=\gamma_t(dx)\,dt.
\]

\vspace{0.2cm}
\noindent
On each cell, the exact discrete energy identity gives
\[
\int_\Omega\widetilde I(u_h(t,x))\,dx
=\mathcal E(t_j)+C_*|\Omega|.
\]
Since $\mathcal E$ is continuous on $[0,T]$, the right-hand side
converges uniformly in time to $\mathcal E(t)+C_*|\Omega|$.
Testing the measure convergence with functions constant in space
therefore gives
\[
\mu_t(\Omega)=\mathcal E(t)+C_*|\Omega|,
\qquad
\int_\Omega\langle\nu,I\rangle\,dx+\gamma_t(\Omega)=\mathcal E(t)
\]
for almost every $t$, proving \eqref{energy-limit}.
By \eqref{energyconsistency}, $U_h$ and the compatible generating
sequence $W_h$ have the same energy-measure limit as $u_h$.
In particular, for every $\varphi\in C_c^1([0,T))$,
\[
\int_{Q_T}\varphi'(t)I(U_h)\,dx\,dt
\longrightarrow
\int_{Q_T}\varphi'(t)\langle\nu,I\rangle\,dx\,dt
+\int_{Q_T}\varphi'(t)\,d\gamma.
\]
Passing to the limit in \eqref{affineenergybalance} yields
\[
E_0\varphi(0)
+\int_{Q_T}\varphi'(t)\langle\nu,I\rangle\,dx\,dt
+\int_{Q_T}\varphi'(t)\,d\gamma
+\int_{Q_T}r\varphi\,dx\,dt=0,
\]
which is \eqref{DMVenergy}. The bound on $\mu_t(\Omega)$ also shows
that $\gamma$ is finite.

}

\vspace{0.4cm}


\vspace{0.4cm}

\noindent\underline{Step 6: Initial entropy trace.} To identify the initial entropy trace, we use
the energy balance together with the support restriction.
For $\eta\geq\eta^0(x)$, assumption (H4) gives
\[
e(F,\eta)
\geq e(F,\eta^0(x))+\theta_0(\eta-\eta^0(x)).
\]
Integrating with respect to $\nu_{t,x}$ and applying
Jensen's inequality to the convex function
$(v,F)\mapsto\frac12|v|^2+e(F,\eta^0(x))$, we obtain
\[
\begin{aligned}
\nuav{I(v,F,\eta)}
\geq
\nuav{\tfrac12|v|^2+e(F,\eta^0(x))}
+\theta_0(\bar\eta-\eta^0)\geq
\frac12|\bar v|^2+e(\bar F,\eta^0(x))
+\theta_0(\bar\eta-\eta^0).
\end{aligned}
\]
The expectations are finite almost everywhere by
\eqref{moments}, the quadratic growth of $e$, and
$\eta^0\in L^2(\Om)$.
Since $\bar\eta\geq\eta^0$ and $\gamma_t\geq0$,
\eqref{energy-limit} implies
\begin{equation}\label{initial-control}
\theta_0\|\bar\eta(t)-\eta^0\|_{L^1(\Om)}
\leq
{\color{black}\mathcal E(t)-\int_\Omega I(\bar v(t,x),\bar F(t,x),\eta^0(x))\,dx},
\end{equation}
for almost every $t$.

\vspace{0.2cm}

{\color{black}\noindent Define
\[
\mathcal I_{\eta^0}(v,F)=
\int_\Omega\left(\frac12|v|^2+e(F,\eta^0(x))\right)\,dx.
\]
This functional is finite and weakly lower semicontinuous on
$L^2(\Omega;\mathbb R^d)\times L^2(\Omega;\mathbb R^{d\times d})$.
The weak continuity of $\bar v$ and
$\bar F$ at zero gives
\[
\liminf_{t\downarrow0}\mathcal I_{\eta^0}(\bar v(t),\bar F(t))
\geq\mathcal I_{\eta^0}(v^0,F^0)=E_0.
\]
Since $\mathcal E(t)\to E_0$, \eqref{initial-control} implies
\[
\operatorname*{ess\,lim}_{t\downarrow0}
\|\bar\eta(t)-\eta^0\|_{L^1(\Omega)}=0.
\]
Choose the nondecreasing representative of $\bar\eta$ from Step 4.
For this representative, the nonnegative quantity
$\|\bar\eta(t)-\eta^0\|_{L^1}$ is nondecreasing, so its essential
right-hand limit is also its ordinary right-hand limit.
Thus \eqref{entropytrace} holds for this representative.} 

\vspace{0.2cm}

It remains to establish~\eqref{weaketa} for test functions reaching $t=0$.
Let $\zeta\in C_c^1([0,T)\times\Omega)$.
Choose a smooth nondecreasing function
$\chi:[0,\infty)\to[0,1]$ such that
$\chi=0$ on $[0,1/2]$ and $\chi=1$ on $[1,\infty)$,
and set $\chi_\varepsilon(t)=\chi(t/\varepsilon)$.
Testing the interior entropy identity with
$\chi_\varepsilon\zeta$ gives
\[
-\int_{Q_T}\bar\eta\,\chi_\varepsilon
\partial_t\zeta\,dx\,dt
-\int_{Q_T}\bar\eta\,\chi_\varepsilon'\zeta\,dx\,dt
=
\int_{Q_T}
r\langle\nu,1/\theta(F,\eta)\rangle
\chi_\varepsilon\zeta\,dx\,dt
+\int_{Q_T}\chi_\varepsilon\zeta\,d\sigma.
\]
By the initial trace~(57), the continuity of $\zeta$, and
$\int_0^T\chi_\varepsilon'(t)\,dt=1$ for sufficiently small
$\varepsilon$, we have
\[
\int_{Q_T}\bar\eta\,\chi_\varepsilon'\zeta\,dx\,dt
\longrightarrow
\int_\Omega\eta^0\zeta(0)\,dx.
\]
Since $\sigma$ is finite on $Q_T$, dominated convergence
applies to the measure term. It also applies to the remaining
terms, since $\bar\eta\in L^1(Q_T)$ and
\[
0\leq r\langle\nu,1/\theta(F,\eta)\rangle
\leq \theta_0^{-1}r\in L^1(Q_T).
\]
Letting $\varepsilon\to0^+$ therefore yields
\[
-\int_{Q_T}\bar\eta\,\partial_t\zeta\,dx\,dt
-\int_\Omega\eta^0\zeta(0)\,dx
=
\int_{Q_T}r\langle\nu,1/\theta(F,\eta)\rangle
\zeta\,dx\,dt
+\int_{Q_T}\zeta\,d\sigma,
\]
which proves~(56). Since $\sigma\geq0$, this also yields
the entropy inequality in Definition~1. This completes
the proof Theorem \ref{thm}. We note this is equivalent to extracting $\sigma$ as a measure on $[0,T]\times \Omega$ and proving that there are no initial atoms, i.e. that $\sigma(\{0\}\times\Omega) = 0$.
\end{proof}

\begin{remark} \rm
The nonnegative measure $\sigma$ in \eqref{weaketa} records the entropy
produced in addition to that supplied by the prescribed heat source. At each time step, the entropy correction converts the mechanical dissipation $D^j$ into internal energy, producing
an entropy increment of $D^j/\theta_\#^j$. The bounds on the corresponding production rates $D^j/(h\theta_\#^j)$ do not exclude concentration in space or time and therefore only converge in the sense of measures to $\sigma$. Its mass on a space-time region quantifies the entropy produced there beyond the contribution of the heat supply. The physical distinction is clearest when $r=0$ where no heat enters the body and the total energy is conserved, yet the entropy may increase. In this sense, $\sigma$ measures irreversibility in the limiting evolution without implying a loss of total energy.
\end{remark}

\vspace{0.2cm}

\noindent  The entropy correction relies on the
pointwise nonnegativity of the convex Bregman
remainder. Replacing convexity by quasiconvexity
in the unchanged construction does not provide
this property: an integrated quasiconvexity
inequality does not imply the pointwise
nonnegativity required by the correction.
The present result is therefore confined to the
convex constitutive class considered here, and
does not establish existence under the general
quasiconvex hypotheses of \cite{GKV2022}, a direction that we are currently investigating.

\vspace{0.2cm}

\noindent The entropy correction and the
incorporation of the prescribed heat supply
distinguish this scheme from the isentropic
construction of \cite{CGT2020}.
The present work provides energy conservation via the corrective step,
and still is consistent with positive entropy production.
The argument proves the mechanical minimization,
the stability estimates, and convergence to a
dissipative measure-valued solution.
It does not establish the smooth-limit convergence
theorem proved for the original scheme in
\cite{CGT2020}; that conclusion would require
a separate relative-energy argument under
suitable additional assumptions.

\section*{Acknowledgements}
The research of KK was partially supported by EPSRC, United Kingdom EP/X038998/1
grant. AV is funded in whole by the Austrian Science Fund (FWF), grant number 10.55776/ESP2735124.

\end{document}